# Representations of Color Lie Superalgebras by Hilbert Series

Shadi Shaqaqha

Yarmouk University, Irbid-Jordan

shadi.s@yu.edu.jo

**Abstract**

The representations of various color Lie superalgebras by Hilbert series are the main topic of this work. The Color Lie superalgebras appear in various branches of mathematics (e.g., topology, algebraic groups, etc.). They are generalized Lie superalgebras. A generating function known as the Hilbert series of color Lie superalgebras which encodes crucial knowledge about the superalgebras representation. In particular, it provides a way to count the number of states in the a given degree. We present a dimension formula that resembles Witt's formula for free color Lie superalgebras, and a specific class of color Lie $p-$superalgebras.

**Keywords;** Hilbert series, color Lie superalgebras, free color Lie superalgebras, restricted color Lie superalgebras, superalgebras representations, vector space.

## 1. INTRODUCTION

Hilbert series of common Lie superalgebra representations is a topic in the field of algebra and representation theory. A Lie superalgebra is a mathematical structure that generalizes the concept of Lie algebra, a vector space equipped with a binary operation that satisfies specific properties. Lie superalgebras are used in various mathematical and physical applications, including quantum mechanics, differential geometry, and string theory. Hilbert series of algebra is a subset of the Hilbert-Poincare series of a graded vector space [1]. Consider $V = \bigoplus_{k=0}^{\infty} V_k$ is a graded vector space such that all subspaces $V_k$ are finite-dimensional. The formal power series in the indeterminate $t$

$$H(V,t) = \sum_{k=0}^{\infty} (\dim V_k) t^k$$

Let $V = U_{k=1}^{\infty} V^k$ be a filtered vector space such that $dim V^k < \infty$ for all $k \in N$. Set $V^0 = 0$. The Hilbert-Poincaré series of V is $H(V) = H(V,t) = \sum_{k=1}^{\infty} \dim \left( V^k / V^{k-1} \right) t^k$. In other words, the Hilbert-Poincaré series for a filtered space V is the same as the associated graded space: $H(V,t) = H(grV,t)$.

Suppose $L = \bigoplus_{n=1}^{\infty} L_n$ be a free Lie algebra of rank $r$. The well-known Witt formula gives the dimensions of homogeneous subspaces $L_n$:

$$dim L_n = \frac{1}{n} \sum_{d|n} \mu(d) r^{\frac{n}{d}},$$

Where $\mu : \mathbb{N} \to \{-1, 0, 1\}$ is the Möbius function defined below. If $n$ is divisible by a prime number's square, we set $\mu(n) = 0$; otherwise, we set $\mu(n) = (-1)^k$, where k is the number of prime divisors of $n$ (with $k = 0$ for $n = 1$, so $\mu(1) = 1$). For homogeneous and multi-homogeneous components of free (color) Lie superalgebras, similar formulas exist. Petrogradsky discovered dimension formulas for free Lie p-algebras [2]. More broadly, suppose $\Lambda$ is a countable abelian semigroup in which every element $\lambda \in \Lambda$ can be written as a sum of other elements only in a finite number of ways. Let $L = \bigoplus_{\lambda \in \Lambda} L_\lambda$ be a $\Lambda$-graded Lie algebra generated freely by $X = \bigcup_{\lambda \in \Lambda} X_\lambda$. Kang and Kim discovered an analog of Witt's formula, known as the character formula, for the dimensions of homogeneous components $L_\lambda, \lambda \in \Lambda$, in [3]. Shaqaqha also worked on free Lie superalgebras and their related formulas [4].

According to Nielsen and Schreier's well-known Theorem, every subgroup of a free group is again free. Shirshov and Witt independently obtained a corresponding result for Lie algebras. On the other hand, subalgebras of the free associative algebra are not always free (for example, $F[x^2, x^3] \subseteq F[x]$ is not free) [5,6].

According to Kukin, $A$ has a filtration (as an algebra) $\bigcup_{i=1}^{\infty} A^i$ if generated by a finitely graded set $X.A^i$, where $A^i$ is spanned by all weight monomials up to $i$. We denote the corresponding series by $H_X(A, t)$ [7]. If $X$ can freely generate $A$, then

$$H_X(A,t) = H(Y,t) = \sum_{i=1}^{\infty} |Y_i| t^i,$$

Where $Y$ denotes the set of all finitely graded monomials in $X$. If $B$ is a subspace of $A$, then the factor-space $A/B$ gains a filtration as well.

$$(A/B)^n = (A^n + B)/B \cong A^n / (B \cap A^n).$$

Petrogradsky defined an operator $\mathcal{E}$ on $Z[[t]]$ (the ring of formal power series in the indeterminate $t$ over $\mathbf{Z}$) as follows in [8]:

$$\varepsilon : \sum_{i=0}^{\infty} a_i t^i \mapsto \prod_{i=0}^{\infty} \frac{1}{(1-t^i)^{a_i}}.$$

Then he presented a formal power series analogue of Schreier's formula for free Lie algebras. He demonstrated that if $L$ is a free Lie algebra generated by a finitely graded set $X$ and $K$ is a subalgebra of $L$, there exists a set of free generators $Z$ of $K$ such that.

$$H(Z) = (H(X) - 1)\,\varepsilon(H(L/K)) + 1.$$

## 2. Representations on Character Formulas for Color Lie Superalgebras

Let $\Lambda$ be a countable additive abelian semigroup such that every element $\lambda \in \Lambda$ can be written as a sum of other elements only in finitely many ways (the finiteness condition). In order to study (color) Lie ($p$-) superalgebras, we fix a homomorphism $\kappa : \Lambda \to \mathbb{Z}_2 = \{\pm 1\}$. This implies that $\Lambda$ can be partitioned as

$$\Lambda = \Lambda_+ \cup \Lambda_-,$$

Where,

$$\Lambda_{\pm} = \{\lambda \in \Lambda \mid \kappa(\lambda) = \pm 1\}.$$

In this section, we consider $\Lambda$-graded color Lie superalgebras $= \bigoplus_{\lambda \in G} L_\lambda$, where for each $g \in G_+$ (respectively, $g \in G_-$), we have $L_g = \bigoplus_{\lambda \in \Lambda_+}(L_\lambda \cap L_g)$ (respectively, $L_g = \bigoplus_{\lambda \in \Lambda_-}(L_\lambda \cap L_g)$). The main purpose in this paper is to derive a dimension formula for the homogeneous subspaces of the free color Lie superalgebras. Also, we will obtain similar results for a certain case of color Lie $p$-superalgebras. $A = \bigoplus_{n \geq 0} A_n$

### 2.1 Characters of Color Lie Superalgebras

Let $U = \bigoplus_{\lambda \in \Lambda} U_\lambda$ be a $\Lambda$-graded space. The character of $U$ is defined by

$$ch_\Lambda U = \sum_{\lambda \in \Lambda} (dim U_\lambda) e^\lambda.$$

It is an element in $\mathbb{Q}[[\Lambda]]$, the completion of the semigroup algebra $\mathbb{Q}[\Lambda]$, whose basis consists of symbols $e^\lambda, \lambda \in \Lambda$ with the multiplication $e^\lambda e^\lambda = e^{\lambda+\mu}$ for all $\lambda, \mu \in \Lambda$. Gradings $U = \bigoplus_{\lambda \in \Lambda} U_\lambda$ and $V = \bigoplus_{\lambda \in \Lambda} V_\lambda$ induce gradings on the spaces $U \oplus V$ and $U \otimes V$:

$$(U \oplus V)_\lambda = U_\lambda \oplus V_\lambda; \quad (U \otimes V)_\lambda = \sum_{\lambda = u + v} (U_u \otimes V_v).$$

By the finiteness condition, the sum above is finite. The following Theorem holds.

**2.1.1. Theorem** $ch_\Lambda(U \oplus V) = ch_\Lambda U + ch_\Lambda V$, and $ch_\Lambda(U \otimes V) = ch_\Lambda U ch_\Lambda V$. A critical special case is $\Lambda = \mathbb{N}$, where $\mathbb{Q}[[\Lambda]]$ is the algebra of formal power series in one variable (without constant term).

*2.2 Characters of Color Lie Superalgebras and Their Enveloping Algebras*

Let $L = L_+ \oplus L_-$ be a free color Lie superalgebra generated by $X$ where $L_\pm = \bigoplus_{\lambda \in \Lambda_\pm} L_\lambda$, with $dimL_\lambda < \infty, \forall \lambda \in \Lambda$ over $F$. The author considered a particular case of our grading in research [9]. The grading by $\Lambda = \Gamma \times G$, where $\Gamma$ is a countable additive abelian semigroup satisfying the following condition: every element $(\alpha, g) \in \Gamma \times G$ can be presented as a sum of other elements only in finitely many ways, and also $\Lambda_+ = \Gamma \times G_+$ and $\Lambda_- = \Gamma \times G_-$. As before, the character of $L$ with respect to $\Lambda$-grading is

$$ch_\Lambda L = \sum_{\lambda \in \Lambda} (dimL_\lambda) e^\lambda, dimL_\pm = \sum_{\lambda \in \Lambda_\pm} (dimL_\lambda) e^\lambda.$$

Note that the universal enveloping algebra is graded by $\bar{\Lambda} = \Lambda \cup \{0\}$. We shall give here the proof of the following formula, established in a study which relates the characters of Lie color superalgebra to that its enveloping algebra [10].

**2.2.1. Lemma** Let $L = L_+ \oplus L_-$ be a $\Lambda$-graded color Lie superalgebra. Then

$$ch_{\bar{\Lambda}} U(L) = \frac{\Pi_{\lambda \in \Lambda_-}(1+ e^\lambda)^{dimL_\lambda}}{\Pi_{\lambda \in \Lambda_+}(1- e^\lambda)^{dimL_\lambda}}.$$

*Proof.* Let $\{e_\lambda \mid \lambda \in \Lambda\}$ be a basis of the positive part $L_+$ and $\{f_\mu \mid \mu \in \Lambda\}$ be a basis of the negative part $L_-$. $U(L)$, as a vector space, is the tensor product of the polynomial algebra $F[\ldots, e_\lambda, \ldots]$ and the Grassmann algebra $\Lambda[\ldots, f_\mu, \ldots]$. Now, the result follows from Theorem 2.1.1.

The super dimension of the homogeneous subspace $L_\lambda$ is defined by

$$sdim L_\lambda = k(\lambda) dim L_\lambda, \quad \lambda \in \Lambda.$$

Note that

$$ch_\Lambda L = \sum_{\lambda \in \Lambda} (sdim L_\lambda) E^\lambda \in \mathbb{Q}[[\Lambda]],$$

where $E^\lambda = \kappa(\lambda) e^\lambda$. It is convenient to define the following operation, called the *twisted dilation*, on $\mathbb{Q}[[\bar\Lambda]]$:

$$[m] : \sum_{\lambda \in \Lambda} f^\lambda E^\lambda \to f^\lambda E^{m\lambda} \quad m \in \mathbb{N}.$$

### 2.2.2. Lemma

1. $f^{[1]} = f$,

2. the dilation $f \mapsto f^{[m]}$ is an endomorphism of the algebra $\mathbb{Q}[[\bar\Lambda]]$,

3. $\left(f^{[m]}\right)^{[n]} = \left(f^{[n]}\right)^{[m]} = f^{[mn]}$ for all $m, n \in \mathbb{N}$.

*Let us define the following two operators over formal series:*

$$\mathcal{E}: \mathbb{Q}[[\Lambda]] \to 1 + \mathbb{Q}[[\Lambda]] : f \mapsto exp\left(\sum_{m=1}^{\infty} \frac{1}{m} f^m\right),$$

$$\mathcal{L}: 1 + \mathbb{Q}[[\Lambda]] \to \mathbb{Q}[[\Lambda]] : f \mapsto \sum_{n=1}^{\infty} \frac{\mu(n)}{n} ln f^n.$$

The following lemma, proved by Petrogradsky in a study, shows that the operators above are similar to the exponential and logarithm [11].

### 2.2.3. Lemma

1. The mappings $\mathcal{E}$ and $\mathcal{L}$ are well-defined and mutually inverse,

2. $\mathcal{E}(f_1 + f_2) = \mathcal{E}(f_1)\mathcal{E}(f_2), f_1, f_2 \in \mathbb{Q}[[\Lambda]]$,

3. $\mathcal{L}(f_1 f_2) = \mathcal{L}(f_1) + \mathcal{L}(f_2), f_1, f_2 \in 1 + \mathbb{Q}[[\Lambda]]$.

Lemma 2.2.1 was used by Petrogradsky to prove the following Theorem.

**2.2.4. Theorem** *Let $L = \bigoplus_{\lambda \in \Lambda} L_\lambda$ be a $\Lambda$-graded color Lie superalgebra, and $U(L)$ be its enveloping algebra [11]. Then*

1. $ch_{\bar{\Lambda}} U(L) = \mathcal{E}(ch_\Lambda L)$,

2. $ch_\Lambda L = \mathcal{L}(ch_{\bar{\Lambda}} U(L))$.

### 2.3 G- Characters of Color Lie Superalgebras and Their Enveloping Algebra

Assume that the $G$-grading on $L$ is determined by the $\Lambda$-grading in the sense that: there exists a homomorphism $\kappa_G : \Lambda \to G$ such that $L_g = \bigoplus_{\substack{\lambda \in \Lambda \\ \kappa_G(\lambda) = g}} L_\lambda$. Define $\upsilon : G \to \mathbb{Z}_2 = \{\pm 1\}$ by $\upsilon(g) = 1$ (Respectively, $-1$) if $g \in G_+$ (respectively, $g \in G_-$). In this case, we can define the $G$-character of $= \bigoplus_{\lambda \in \Lambda} L_\lambda$, where $dim L_\lambda < \infty$ for all $\lambda \in \Lambda$, as follows

$$ch_\Lambda L = \sum_{\lambda \in \Lambda} (dim L_\lambda) K_G(\lambda) e^\lambda \in \mathbb{Q}[G][[\Lambda]],$$

where $\mathbb{Q}[G]$ is the group algebra of $G$ with coefficients in $\mathbb{Q}$ and $\mathbb{Q}[G][[\Lambda]]$ is the completion of the semigroup algebra $\mathbb{Q}[G][\Lambda]$. For $\lambda \in \Lambda$, we set $sdim L_\lambda = \upsilon(\kappa_G(\lambda)) dim L_\lambda$ and color super dimension $csdim L_\lambda = \kappa_G(\lambda) sdim L_\lambda$. Now, the twisted dilation is defined by

$$[m]: \sum_{\lambda \in \overline{\Lambda}} r_\lambda g_\lambda E^\lambda \to \sum_{\lambda \in \overline{\Lambda}} r_\lambda g_{\lambda^m} E^{m\lambda} \ r_\lambda \in \mathbb{Q}, \lambda_{g_\lambda} \in G, and \ m \in \mathbb{N},$$

Where $E_\lambda = v(\kappa_G(\lambda))\kappa_G(\lambda)e^\lambda$. The character of $L$ can also be written as

$$ch_\Lambda L = \sum_{\lambda \in \Lambda} (sdimL_\lambda) E^\lambda.$$

We have the following properties of the twisted dilation operator.

### 2.3.1. Lemma

1. The dilation $f \mapsto f^{[m]}$ is an endomorphism of the algebra $Q[G][[\Lambda]]$,

2. $\left(f^{[m]}\right)^{[n]} = \left(f^{[n]}\right)^{[m]} = f^{[mn]}$ for all $m, n \in \mathbb{N}$.

3. $\left(\sum_{\lambda \in \Lambda} r_\lambda g_\lambda e^\lambda\right)^{[m]} = \sum_{\lambda \in \Lambda} r_\lambda g_\lambda^m \left(v(\kappa_G(\lambda))\right)^{m+1} e^{m\lambda}, r_\lambda \in \mathbb{Q}, g_\lambda \in G.$

*Proof.* It is clear that the first two properties hold. Hence it remains to prove the last claim.

$$\left(\sum_{\lambda \in \Lambda} r_\lambda g_\lambda e^\lambda\right)^{[m]} = \left(\sum_{\lambda \in \Lambda} r_\lambda g_\lambda^v \left(\kappa_G(\lambda)\right)\left(\kappa_G(\lambda)\right)^{-1} E^\lambda\right)^{[m]}$$

$$= \sum_{\lambda \in \Lambda} r_\lambda v\left(\kappa_G(\lambda)\right) g_\lambda^m \left(\kappa_G(\lambda)\right)^{-m} E^{m\lambda}$$

$$= \sum_{\lambda \in \Lambda} r_\lambda v\left(\kappa_G(\lambda)\right) g_\lambda^m \left(\kappa_G(\lambda)\right)^{-m} v\left(\kappa_G(m\lambda)\right) \kappa_G(m\lambda) e^{m\lambda}$$

$$= \sum_{\lambda \in \Lambda} r_\lambda (v(\kappa_G(\lambda)))^{m+1} g_\lambda^m e^m.$$

We introduce the following two operators over formal power series:

$$\mathcal{E}_G: \mathbb{Q}[G][[\Lambda]] \to 1 + \mathbb{Q}[G][[\Lambda]] : f \mapsto \exp\left(\sum_{m=1}^{\infty} \frac{1}{m} f^{[m]}\right),$$

$$\mathcal{L}_G : 1 + \mathbb{Q}[G][[\Lambda]] \to \mathbb{Q}[G][[\Lambda]] : f \mapsto \sum_{n=1}^{\infty} \frac{\mu(n)}{n} \ln f^{[n]}.$$

We can easily prove the following lemma.

### 2.3.2. Lemma

1. The mappings $\mathcal{E}_G$ and $\mathcal{L}_G$ are well-defined and are mutually inverse.

2. $\mathcal{E}_G(f_1 + f_2) = \mathcal{E}_G(f_1)\,\mathcal{E}_G(f_2), f_1, f_2 \in \mathbb{Q}[G][[\Lambda]]$,

3. $\mathcal{L}_G(f_1 f_2) = \mathcal{L}_G(f_1) + \mathcal{L}_G(f_2), f_1, f_2 \in 1 + \mathbb{Q}[G][[\Lambda]]$.

**2.3.3. Theorem** *Let $L = \bigoplus_{\lambda \in \Lambda} L_\lambda$ be a $\Lambda$-graded color Lie superalgebra and $U(L)$ be its enveloping algebra. Then*

1. $ch_{\overline{\Lambda}}^G U(L) = \mathcal{E}_G(ch_{\Lambda}^G L)$,

2. $ch_{\Lambda}^G L = \mathcal{L}_G(ch_{\overline{\Lambda}}^G U(L))$.

*Proof.* According to PBW-Theorem, we have

$$ch_{\overline{\Lambda}}^G U(L) = \prod_{\lambda \in \Lambda} (1 - E^\lambda)^{-sdimL_\lambda}.$$

Then we see that

$$ch_{\overline{\Lambda}}^G U(L) = \exp\left(-\sum_{\lambda \in \Lambda} (sdimL_\lambda)(1 - E^\lambda)\right).$$

Using $\ln(1+x) = \sum_{n=1}^{\infty} (-1)^{n+1} \frac{x^n}{n}$, we obtain

$$ch^G_{\underline{\Lambda}} U(L) = \exp\left(-\sum_{\lambda \in \Lambda} (sdimL_\lambda) \sum_{m=1}^{\infty} \frac{E^{m\lambda}}{m}\right).$$

Then,

$$ch^G_{\underline{\Lambda}} U(L) = \exp\left(\sum_{m=1}^{\infty} \frac{1}{m} \sum_{\lambda \in \Lambda} (sdimL_\lambda) E^{m\lambda}\right)$$

$$= \exp\left(\sum_{m=1}^{\infty} \frac{1}{m} \left(ch^G\right)^{[m]} L\right)$$

$$= \mathcal{E}_G ch^G_{\underline{\Lambda}} L.$$

To prove the second relation, note that

$$ch^G_{\underline{\Lambda}} L = \mathcal{L}_G \mathcal{E}_G \left(ch^G_{\underline{\Lambda}} L\right) = \mathcal{L}_G (\mathcal{E}_G \left(ch^G_{\underline{\Lambda}} L\right)) = \mathcal{L}(ch^G_{\underline{\Lambda}} U(L)).$$

## 2.4 Character Formula of Free Color Lie Superalgebras

By a $\Lambda$-graded set, we mean a disjoint union $X = \bigcup_{\lambda \in \Lambda} X_\lambda$. If in addition, we have $|X_\lambda| < \infty$ for all $\lambda \in \Lambda$, then we define its character

$$ch_\Lambda X = \sum_{\lambda \in \Lambda} |X_\lambda| e^\lambda \in \mathbb{Q}[[\Lambda]],$$

For an element $x \in X_\lambda \subseteq X$, we say $\Lambda$-weight of $x$ is $\lambda$, and we write $wt_\Lambda x = \lambda$. We call such a set $\Lambda$-*finitely graded* (if $\Lambda = \mathbb{N}$, then we say $X$ is a finitely graded set). For any monomial $y = x_1 \ldots x_n$, where $x_j \in X$, we set $wt_\Lambda y = wt_\Lambda x_1 + \ldots + wt_\Lambda x_n$. Suppose $Y$ is a set of all monomials (associative, Lie, ...) in $X$. We denote

$$Y_\lambda = \{y \in Y \,|\, wt_\Lambda y = \lambda\}.$$

Also, the $\Lambda$-generating function of $Y$ is

$$ch_\Lambda Y = \sum_{\lambda \in \Lambda} |Y_\lambda| e^\lambda \in \mathbb{Q}[[\Lambda]],$$

**2.4.1. Lemma** Let $X = \bigcup_{\lambda \in \Lambda} X_\lambda$ be a $\Lambda$-graded set with $|X_\lambda| < \infty$, $\lambda \in \Lambda$, and let $F(X)$ be the free associative algebra generated by $X$. Then

$$ch_{\bar{\Lambda}} F(X) = \sum_{n=0}^{\infty} (ch_\Lambda X)^n = \frac{1}{1 - ch_\Lambda X}.$$

Petrogradsky proved the following Theorem in the context of Lie superalgebras in [11].

**2.4.2. Theorem** Let $L = L(X)$ be the free color Lie superalgebra generated by a $\Lambda$-graded set $X = \bigcup_{\lambda \in \Lambda} X_\lambda$ with $|X_\lambda| < \infty$ for all $\lambda \in \Lambda$. Then

$$ch_\Lambda L(X) = -\sum_{n=1}^{\infty} \frac{\mu(n)}{n} \ln(1 - ch_\Lambda^{[N]} X).$$

*Proof.* The universal enveloping algebra $U(L)$ is isomorphic to the free associative algebra $F(X)$ generated by $X$. Thus

$$ch_{\bar{\Lambda}} U(L) = \frac{1}{1 - ch_\Lambda X}.$$

Applying Theorem 2.2.4, we have

$$ch_\Lambda L = \mathcal{L}\big(\mathcal{E}(ch_\Lambda L)\big) = \mathcal{L}\left(\frac{1}{1 - ch_\Lambda X}\right) = -\sum_{n=1}^{\infty} \frac{\mu(n)}{n} \ln(1 - ch_\Lambda^{[N]} X),$$

as desired.

We are going to discuss several corollaries of the above result.

If $|G| = r$, we can make any finite set $X$ a $\Lambda$-graded set for $\bar{\Lambda} = \mathbb{N}_0^r$. Write $G = G_+ \cup G_-$ where $G_+ = \{g_1, \ldots, g_k\}$ and $G_- = \{g_{k+1}, \ldots, g_r\}$ (of course, $|G| = |G_+|$ or $|G_+| = |G_-|$) is an abelian group, and $L$ is a free color Lie superalgebra freely generated by a set $X = X_{g1} \cup \cdots \cup X_{gr}$, with $|X_{gi}| = s_i \geq 1, i = 1,\ldots,r$. Consider the case $\Lambda = \mathbb{N}_0^r$. We define a weight function

$$wt: X \to \mathbb{N}_0^r : x \mapsto \lambda_i, \text{ for } i = 1,\ldots,r \text{ and } x \in X_{gi},$$

where $\lambda_i = (0,\ldots,0,1,0,\ldots,0)$ with 1 in the ith place. We define the homomorphism $\kappa: \mathbb{N}_0^r \to \mathbb{Z}_2 = \{\pm 1\}$ by $\kappa(\lambda_i) = 1$ for $1 \leq i \leq k$ and $\kappa(\lambda_i) = -1$ for $k+1 \leq i \leq r$. We denote $t_i = e^{\lambda_i}$, so the algebra $\mathbb{Q}[[\overline{\Lambda}]]$ turns into the formal power series ring $\mathbb{Q}[[t]] = \mathbb{Q}[[t_1,\ldots,t_r]]$. In this case, the character of a $\Lambda$-graded Lie superalgebra, $L$, is the multivariable Hilbert-Poincaré series, $H(L,t) = H(L; t_1,\ldots,t_r)$, of $L$. We have the following result.

**2.4.3. Corollary** Suppose *that $G = G_+ \cup G_-$ is an abelian group, where $G_+ = \{g_1,\ldots,g_k\}$ and $G_- = \{g_{k+1},\ldots,g_r\}$ ($r = k$ or $r = 2k$), and $L$ is a free color Lie superalgebra freely generated by a set $X = X_{g1} \cup \cdots \cup X_{gr}$ with $|X_{gi}| = si \geq 1, i = 1,\ldots,r$. Then*

$$H(L; t_1\ldots,t_k, t_{k+1},\ldots,t_r) = -\sum_{n=1}^{\infty} \frac{\mu(n)}{n} \ln\left(1 - \sum_{i=1}^{k} s_i t_i^n + \sum_{j=k+1}^{r} s_j(-t_j)^n\right).$$

*Proof.* In this case $ch_\Lambda L = \sum_{i=1}^{r} s_i t_i$, and so $ch^{[n]}X = \sum_{i=1}^{k} s_i t_i^n - \sum_{j=k+1}^{j} s_j(-t_j)^n$.

The formula follows from Theorem 2.4.2

The weight function $wt: X \to \mathbb{N}_0^r$ defines the multidegree $\alpha = (\alpha_1,\ldots,\alpha_r) \in \mathbb{N}_0^r$ for elements of $L$, and the degree $|\alpha| = \alpha_1 + \cdots + \alpha_r$. Also, we write $|\alpha|_+ = \alpha_1 + \cdots + \alpha_k$ and $|\alpha|_- = \alpha_{k+1} + \cdots + \alpha_r$. By $n|\alpha$ we denote that n divides all components $\alpha_i$ of $\alpha$. Then we have the following result.

**2.4.4. Corollary** Suppose $G = G_+ \cup G_-$ and $L = L(X)$ as in Corollary 2.4.3. Then

$$\dim L_\alpha = \frac{(-1)^{|\alpha|_-}}{|\alpha|} \sum_{n|\alpha} \mu(n) \frac{\left(\frac{|\alpha|}{n}\right)! \, (-1)^{\frac{|\alpha|_-}{n}}}{\left(\frac{\alpha_1}{n}\right)! \cdots \left(\frac{\alpha_r}{n}\right)!} s_1^{\frac{\alpha_1}{n}} \cdots s_r^{\frac{\alpha_r}{n}}.$$

*In particular, if L is a free Lie algebra, we get the classical Witt's formula.*

*Proof.* We apply the formula for $H(L; t_1, \ldots, t_r)$ from the corollary above. We have

$$H(L; t) = -\sum_{n=1}^{\infty} \frac{\mu(n)}{n} \ln\left(1 - \sum_{i=1}^{k} s_i t_i^n + \sum_{j=k+1}^{r} s_j(-t_j)^n\right)$$

$$= \sum_{n=1}^{\infty} \frac{\mu(n)}{n} \sum_{s=1}^{\infty} \frac{(s_1 t_1^n + \cdots + s_k t_k^n - s_{k+1}(-t_{k+1})^n)^s}{s}.$$

Applying the multinomial formula, we get

$$H(L; t) = \sum_{n=1}^{\infty} \frac{\mu(n)}{n} \sum_{s=1}^{\infty} \frac{1}{s} \sum_{|\beta|=s} \frac{|\beta|!}{\beta_1! \cdots \beta_r!} (s_1 t_1^n)^{\beta_1} \cdots (s_k t_k^n)^{\beta_k} ((-s_{k+1})(-t_{k+1}^n))^{\beta_{k+1}}$$

$$\cdots ((-s^r)(-t_r^n))^{\beta_r}.$$

Hence,

$$H(L; t) = \sum_{n=1}^{\infty} \frac{\mu(n)}{n} \sum_{s=1}^{\infty} \frac{1}{s} \sum_{|\beta|=s} \frac{|\beta|!(-1)^{(n+1)|\beta|_-}}{\beta_1! \cdots \beta_r!} s_1^{\beta_1} \cdots s_r^{\beta_r} t_1^{n\beta_1} \cdots t_r^{n\beta_r}$$

$$= \sum_{\alpha \in \mathbb{N}_0^r \setminus \{0\}} \frac{1}{|\alpha|} \sum_{n|\alpha} \mu(n) \frac{\left(\frac{|\alpha|}{n}\right)!(-1)^{|\alpha|_-+\frac{|\alpha|_-}{n}}}{\left(\frac{\alpha_1}{n}\right)! \cdots \left(\frac{\alpha_r}{n}\right)!} s_1^{\frac{\alpha_1}{n}} \cdots s_r^{\frac{\alpha_r}{n}} t_1^{\alpha_1} \cdots t_r^{\alpha_r}.$$

On the other hand, $(L; t) = \sum_{\alpha \in \mathbb{N}_0^r \setminus \{0\}} \dim L_\alpha t^\alpha$. Therefore

$$\dim L_\alpha = \frac{(-1)^{|\alpha|_-}}{|\alpha|} \sum_{n|\alpha} \mu(n) \frac{\left(\frac{|\alpha|}{n}\right)!(-1)^{\frac{|\alpha|_-}{n}}}{\left(\frac{\alpha_1}{n}\right)! \cdots \left(\frac{\alpha_r}{n}\right)!} s_1^{\frac{\alpha_1}{n}} \cdots s_r^{\frac{\alpha_r}{n}},$$

as desired.

Let *X* be a finite generating set of the free color Lie superalgebra *L(X)* with the weight functions

$$wt : X \to \mathbb{N}^2,$$

defined by

$$x \mapsto (1,0) \text{ if } x \in X_+ \text{ and } x \mapsto (0,1) \text{ if } x \in X_-.$$

If we denote $t_+ = e^{(1,0)}$ and $t_- = e^{(0,1)}$, then the algebra $\mathbb{Q}[[\mathbb{N}_0^2]]$ is the formal power series ring $\mathbb{Q}[[t_+, t_-]]$. We have the following corollary.

**2.4.5. Corollary** Let $L = L(X)$ be a free color Lie superalgebra freely generated by the set $X = X_+ \cup X_-$, where $X_+ = \{x_1, \ldots, x_k\}$ and $X_- = \{x_{k+1}, \ldots, x_r\}$. Then

1. $H(L; t_+, t_-) = -\sum_{n=1}^{\infty} \frac{\mu(n)}{n} \ln(1 - kt_+^n + (r-k)(-t_-)^n).$

2. $H(L, t) = H(L; t_+, t_-)|_{t_+ = t_- = t} = -\sum_{n=1}^{\infty} \frac{\mu(n)}{n} \ln(1 - (k - (-1)^n (r-k))t^n).$

**2.4.6. Corollary** $L(X)$ be a free color Lie superalgebra freely generated by the set $X = X_+ \cup X_-$, where $X_+ = \{x_1, \ldots, x_k\}$ and $X_- = \{x_{k+1}, \ldots, x_r\}$. Consider the weight function $wt : X \to \mathbb{N}; x \mapsto 1$. Then

$$dim L_n = \frac{1}{n} \sum_{m|n} \mu(m)\left(k - (-1)^m(r-k)\right)^{\frac{n}{m}}.$$

Let us return to the general setting. Let $\Lambda$ and $\Gamma$ be two additive abelian semigroups satisfying the finiteness condition, $\kappa : \Lambda \to \mathbb{Z}_2$ and $\kappa' : \Gamma \to \mathbb{Z}_2$ are homomorphisms. Suppose that $\varphi : \Lambda \to \Gamma$ is a semigroup homomorphism such that $\kappa = \kappa' \circ \varnothing$ and for each $\gamma \in \Gamma$ the set $\{\lambda \in \Lambda \mid \varnothing(\lambda) = \gamma\}$ is finite. Let $L = L_+ \oplus L_-$ be a free $\Lambda$-graded algebra generated by $X = \bigcup_{\lambda \in \Lambda} X_\lambda$. Using the homomorphism $\varnothing$, we can also regard $L$ as $\Gamma$-graded. Then

$$ch_\Gamma L = \sum_{\gamma \in \Gamma} dim L_\gamma e^\gamma = \left( \sum_{\substack{\lambda \in \Lambda \\ \varphi(\lambda) = \gamma}}^{j} dim L_\lambda \right) e^\lambda \cdots (2.1).$$

Now, we consider the case where $\Lambda = \mathbb{N} \times G$. Such a situation can be obtained from the grading given in Theorem 2.4.3 by taking the grading.

$$wt : X \to \mathbb{N} \times G : x \mapsto (1, g_i) \text{ for } x \in X_{g_i}.$$

For such grading, we will use superscripts instead of subscripts. As a result, we have the following corollary.

**2.4.7. Corollary** $dimL^{(n,g)} = \sum_{\substack{\alpha_1+\cdots\alpha_r=n \\ g_1^{\alpha_1}\cdots g_r^{\alpha_r}=g}} dimL_{(\alpha_1,\ldots,\alpha r)}$.

*Proof.* The result is the formula 2.1 applied to

$$\varphi: \mathbb{N}^r \to \mathbb{N} \times G: \lambda_i \mapsto (1, g_i).$$

**2.4.8. Example** Consider the free $(\mathbb{Z}_2 \oplus \mathbb{Z}_2, \gamma)$−color Lie superalgebra $L = L(X)$ over the field $F = \mathbb{C}$ where

$$\gamma: (\mathbb{Z}_2 \oplus \mathbb{Z}_2) \times (\mathbb{Z}_2 \oplus \mathbb{Z}_2) \to \mathbb{C}^*: ((a_1, a_2), (b_1, b_2)) \mapsto (-1)^{(a_1+a_2)(b_1+b_2)}.$$

Hence, $G_+ = \{(0,0), (1,1)\}$ and $G_- = \{(0,1), (1,0)\}$. Let $g_1 = (0,0), g_2 = (1,1), g_3 = (0,1)$, and $g_4 = (1,0)$, and let $|X_{g_1}| = 1, |X_{g_2}| = 2, |X_{g_3}| = |X_{g_4}| = 1$. According to the above Theorem, we have

$$dimL^{(3,(1,1))} = dimL_{(0,3,0,0)} + dimL_{(2,1,0,0)} + dimL_{(1,0,1,1)}.$$

Now, if we apply the formula given in Corollary 2.4.4, we have

$$dimL_{(0,3,0,0)} = \frac{(-1)^0}{3}\left(u(1)\frac{(3!)(-1)^0}{3!}2^3 + u(3)\frac{(1!)(-1)^0}{1!}2^1\right) = 2.$$

Similarly, we obtain $dimL_{(2,1,0,0)} = 2$, and $dimL_{(1,0,1,1)} = 2$. Hence $dimL^{(3,(1,1))} = 2 + 2 + 2 = 6$.

## 2.5 Characters of Free Restricted Color Lie Superalgebras

Let $L = L_+ \oplus L_-$ be a free color restricted Lie superalgebra generated by $X$ where $L_\pm = \oplus_{\lambda \in \Lambda_\pm} L_\lambda$ with $\dim L_\lambda < \infty \; \forall \lambda \in \Lambda$ over a field $F$. We can now deduce the formula that relates the character of Lie color $p$-superalgebra to that of its restricted enveloping algebra.

**2.5.1. Lemma** *Let $L = L_+ \oplus L_-$ be a $\Lambda$-graded color Lie p-superalgebra. Then*

$$ch_{\bar{\lambda}} u(L) = \Pi_{\lambda \in \Lambda_-} (1 + e^\lambda)^{\dim L_\lambda} \Pi_{\lambda \in \Lambda_+} (1 + \cdots + e^{(p-1)\lambda})^{\dim L_\lambda}.$$

*Proof.* For color Lie p-superalgebras, the PBW-theorem must be used, as in Lemma 2.2.1. The specifics are omitted.

The remainder of this section will look at a -graded color Lie p-superalgebra satisfying $G = G_+$; remember that the ordinary restricted Lie algebra is a special case. (Recall that color Lie p-superalgebras are also known as the color Lie p-algebras.)

Petrogradsky has defined functions $1_p, \mu_p : \mathbb{N} \to \mathbb{N}$ by:

$$1_p(n) = \begin{cases} 1, & \text{if } (p,n) = 1 \\ 1 - p & \text{if } (p,n) = p, \end{cases}$$

and

$$\mu_p(n) = \begin{cases} \mu(n), & \text{if } (p,n) = 1 \\ \mu(m)(p^s - p^{s-1}), & \text{if } n = mp^s, (p,m) = 1, s \geq 1. \end{cases}$$

Recall that a function $f : \mathbb{N} \to \mathbb{N}$ is multiplicative if $f(nm) = f(n)f(m)$ for any coprime $n, m$. One can easily show that $1_p$ and $\mu_p$ are multiplicative functions. In addition, we have the following property [12].

**2.5.2 Lemma** $\sum_{ab=n} 1_p(b) \mu_p(a) = 0$ for all $n > 1$.

Proof. We fill in the details of the proof in a study. First, we assume n is not divisible by p. Let a, b ∈ ℕ with $ab = n$. Then $a$ and $b$ are not divisible by $p$. Hence $1_p(b) = 1$ and $\mu_p(a) = \mu(a)$ [12]. Now, the statement follows from the property of the Möbius function. Next, we suppose n is divisible by p. Write $n = n'p^k$, $k \geq 1$, where $n'$ is not divisible by $p$. For all $a, b \in \mathbb{N}$ with $ab = n$, we write accordingly $a = a'p^r$ and $b = b'p^s$ where $r + s = k$. Then

$$\sum_{ab=n} 1_p(b)\mu_p(a) = \sum_{a'b'=n'} 1_p(b')\mu_p(a') \sum_{r+s=k} 1_p(p')\mu_p(p')$$

$$= \sum_{a'b'=n'} \mu(a')(1(p^k - p^{k-1}) + (1-p)(p^{k-1} - p^{k-2}) + \ldots + (1-p)1)$$

$$= \sum_{a'b'=n'} \mu s(a')((p^k - p^{k-1}) + (1-p)(p^{k-1} - 1)$$

$$+ (1-p))$$

$$= 0,$$

where in the first line, we used the fact that $1_p$ and $u_p$ are multiplicative functions.

We introduce the following two operators on formal series, which were defined in the case of $\bar{\Lambda} = \mathbb{N}_0^m$.

$$\mathcal{E}_p : \mathbb{Q}[[\Lambda]] \rightarrow 1 + \mathbb{Q}[[\Lambda]] : f \mapsto \exp\left(\sum_{m=1}^{\infty} \frac{1_{p\,(m)}}{m} f^{[m]}\right),$$

$$\mathcal{L}_p : 1 + \mathbb{Q}[[\Lambda]] \rightarrow Q[[\Lambda]] : f \mapsto \sum_{n=1}^{\infty} \frac{\mu_{p\,(n)}}{n} \ln f^{[n]}.$$

Now we show that these operators are similar to the exponential and logarithm.

**2.5.3. Theorem** *Let $L = L(X)$ be the free color Lie p-algebra ($G = G_+$) generated by a $\Lambda$-graded set $X = \cup_{\lambda \in \Lambda} X_\lambda$, with $|X_\lambda| < \infty$ for all $\lambda \in \Lambda = \Lambda_+$. Then*

$$ch_\Lambda L(X) = -\sum_{n=1}^{\infty} \frac{u_p(n)}{n} \ln\left(1 - ch_{[n]}X\right).$$

**2.5.4. Lemma**

1. The maps $\mathcal{E}_p$ and $\mathcal{L}_p$ are well defined and mutually inverse,

2. $\mathcal{E}_p(f_1 + f_2) = \mathcal{E}_p(f_1)\mathcal{E}_p(f_2), f_1, f_2 \in \mathbb{Q}[[\Lambda]],$

3. $\mathcal{L}_p(f_1 f_2) = \mathcal{L}_p(f1) + \mathcal{L}_p(f_2), f_1, f_2 \in 1 + \mathbb{Q}[[\Lambda]].$

*Proof.* It follows from the finiteness condition of $\Lambda$ that $\mathcal{E}_p$ and $\mathcal{L}_p$ are well defined. Let $f \in \mathbb{Q}[[\Lambda]]$. Then

$$\mathcal{L}_p(\mathcal{E}_p(f)) = \mathcal{L}_p\left(\exp\left(\sum_{m=1}^{\infty} \frac{1_p(m)}{m} f^{[m]}\right)\right) \text{ (Definition of } \mathcal{E}_p)$$

$$= \sum_{n=1}^{\infty} \frac{\mu_p(n)}{n} \ln\left(\exp\left(\sum_{m=1}^{\infty} \frac{1_p(m)}{m} f^{[m]}\right)\right)^{[n]} \text{ (Definition of } \mathcal{L}_p)$$

$$= \sum_{n=1}^{\infty} \frac{\mu_p(n)}{n} \ln\left(\prod_{m=1}^{\infty} \exp\left(\frac{1_p(m)}{m} f^{[m]}\right)\right)^{[n]}$$

$$= \sum_{n=1}^{\infty} \frac{\mu_p(n)}{n} \ln\left(\prod_{m=1}^{\infty} \exp\left(\frac{1_p(m)}{m} f^{[m]}\right)^{[n]}\right) \text{ (Lemma 2.2.2)}$$

$$= \sum_{n=1}^{\infty} \frac{\mu_p(n)}{n} \sum_{m=1}^{\infty} \frac{1_p(m)}{m} (f^{[m]})^{[n]}$$

$$= \sum_{n-1}^{\infty} \sum_{m-1}^{\infty} \frac{f^{[mn]}}{mn} 1_p(m)\mu_p(n) \text{ (Lemma 2.2.2)}$$

$$= \sum_{k=1}^{\infty} \frac{f^{[k]}}{k} \sum_{mn=k} 1_p(m)\mu_p(n)$$

$$= f^{[1]} \text{ (Lemma 2.5.2)}$$

$$= f.$$

In a similar way, we can prove $\mathcal{E}_p(\mathcal{L}_p(f)) = f, f \in 1 + \mathbb{Q}[[\Lambda]]$. The relations (2) and (3) are clear.

**2.5.5. Theorem** *Let $L = \bigoplus_{\lambda \in \Lambda} L_\lambda$ be a $\Lambda$-graded color Lie p-algebra (G = G+) and u(L) be its restricted enveloping algebra. Then*

1. $ch_{\bar{\Lambda}} u(L) = \mathcal{E}_p(ch_\Lambda L),$

2. $ch_\Lambda u(L) = \mathcal{L}_p(ch_{\bar{\Lambda}} u(L)).$

*Proof.*

1. By Lemma 2.5.1, we have

$$ch_{\bar{\Lambda}} u(L) = \prod_{\lambda \in \Lambda} (1 + e^\lambda + \cdots + e^{(p-1)\lambda})^{dim L^\lambda}.$$

Now, as $(1 - e^{p\lambda}) = (1 - e^{\lambda})(1 + e^{\lambda} + \cdots + e^{(p-1)\lambda})$, $ch_{\bar{\Lambda}}u(L)$ can be written as:

$$ch_{\bar{\Lambda}}u(L) = \prod_{\lambda \in \Lambda}\left(\frac{1-e^{p\lambda}}{1-e^{\lambda}}\right)^{dimL^{\lambda}}.$$

Therefore,

$$ch_{\bar{\Lambda}}u(L) = \exp\left(\sum_{\lambda \in \Lambda} dimL^{\lambda}\left((-\ln(1-e^{\lambda}) + \ln(1-e^{p\lambda}))\right)\right).$$

Using $\ln(1+x) = \sum_{n=1}^{\infty}(-1)^{n+1}\frac{x^n}{n}$, we obtain

$$ch_{\bar{\Lambda}}u(L) = \exp\left(\sum_{\lambda \in \Lambda} dimL_{\lambda}\left(\sum_{n=1}^{\infty}\frac{e^{n\lambda}}{n} - \sum_{n=1}^{\infty}\frac{e^{pn\lambda}}{n}\right)\right).$$

Then we see that

$$\begin{aligned}
ch_{\bar{\Lambda}}u(L) &= \exp\left(\sum_{\lambda \in \Lambda} dimL_{\lambda}\left(\sum_{n=1}^{\infty}\frac{e^{n\lambda}}{n} - \sum_{n=1}^{\infty}\frac{e^{pn\lambda}}{n}\right)\right) \\
&= \exp\left(\sum_{\lambda \in \Lambda} dimL_{\lambda}\left(\sum_{n=1,p|n}^{\infty}\frac{e^{n\lambda}}{n} + \sum_{n=1}^{\infty}\left(\frac{e^{np\lambda}}{np} - \frac{e^{np\lambda}}{n}\right)\right)\right). \\
&= \exp\left(\sum_{\lambda \in \Lambda} dimL_{\lambda}\left(\sum_{n=1,p|n}^{\infty}\frac{e^{n\lambda}}{n} + \sum_{n=1}^{\infty}\left(\frac{e^{np\lambda} - pe^{np\lambda}}{np}\right)\right)\right) \\
&= \exp\left(\sum_{\lambda \in \Lambda} dimL_{\lambda}\left(\sum_{n=1}^{\infty}e^{n\lambda}\frac{1_p(n)}{n}\right)\right) \\
&= \exp\left(\sum_{n=1}^{\infty}\frac{1_p(n)}{n}\sum_{\lambda \in \Lambda} dimL_{\lambda}\, e^{n\lambda}\right) \\
&= \exp\left(\sum_{n=1}^{\infty}\frac{1_p(n)}{n}(ch_{\Lambda}L)^{[n]}\right) \\
&= \mathcal{E}_p(ch_{\Lambda}L).
\end{aligned}$$

2. This relation follows directly from Lemma 2.5.4 and (1):

$$ch_{\Lambda}L = \mathcal{L}_p\mathcal{E}_p(ch_{\Lambda}L) = \mathcal{L}_p(\mathcal{E}_p(ch_{\Lambda}L)) = \mathcal{L}_p(ch_{\bar{\Lambda}}u(L)).$$

**2.5.6. Remark** One can also extend the definition of $\mathcal{E}_p$ to the general case $\Lambda = \Lambda_+ \cup \Lambda_-$ as follows:

$$\mathcal{E}_p: \mathbb{Q}[[\Lambda]] \to 1 + \mathbb{Q}[[\Lambda]]: f = f_+ + f_- \mapsto \exp\left(\sum_{m=1}^{\infty}\frac{1_p(m)}{m}f_+^{[m]}\right)\exp\left(\sum_{n=1}^{\infty}\frac{1}{n}f_-^{[n]}\right).$$

Again $\mathcal{E}_p$ is well defined operator. Also, it is easy to see that

1. $\mathcal{E}_p(f_1 + f_2) = \mathcal{E}_p(f_1) \mathcal{E}_p(f_2), f_1, f_2 \in \mathbb{Q}[[\Lambda]]$,

2. $ch_{\bar{\Lambda}} u(L) = \mathcal{E}_p(ch_\Lambda L)$.

**2.5.7. Theorem** *Let $L = L(X)$ be the free color Lie p-algebra ($G = G_+$) generated by a $\Lambda$-graded set $X = \cup_{\lambda \in \Lambda} X_\lambda$, with $|X_\lambda| < \infty$ for all $\lambda \in \Lambda = \Lambda_+$. Then*

$$ch_\Lambda L(X) = -\sum \frac{\mu_p(n)}{n} \ln\left(1 - ch_\Lambda^{[n]} X\right)_{n=1}^\infty.$$

*Proof.* For the restricted color Lie superalgebra $L = L(X)$, we denote the restricted enveloping algebra of $L$ by $u(L)$. Let $F(X)$ be the free associative algebra on X. It is well known that $u(L(X))$ is isomorphic to $F(X)$[13]. Thus,

$$ch_{\bar{\Lambda}} u(L) = \frac{1}{(1 - ch_\Lambda X)}.$$

Using Theorem 2.5.5, we get

$$ch_\Lambda L = \mathcal{L}_p ch_{\bar{\Lambda}} u(L) = \mathcal{L}_p \frac{1}{(1 - ch_\Lambda X)} = -\sum \frac{\mu_p(n)}{n} \ln\left(1 - ch_\Lambda^{[n]} X\right)_{n=1}^\infty.$$

**2.5.8. Corollary** *Let $L = L(X)$ be the free color Lie p-algebra generated by at most countable set $X = \{x_i \mid i \in I\}$. Then*

$$H(L, ti \mid i \in I) = -\sum_{n=0}^\infty \frac{\mu_p(n)}{n} \ln\left(1 - \sum_{i \in I} t_i^n\right).$$

*In particular, if L is generated by $X = \{x_1, \ldots, x_r\}$, then*

$$H(L; t_1, \ldots, t_r) = -\sum_{n=1}^\infty \frac{\mu_p(n)}{n} \ln(1 - t_1^n - \cdots - t_r^n).$$

Consider the particular case $\Lambda = \mathbb{N}$ and $wt : X \to \mathbb{N} : x \mapsto 1$. Then we have the following result.

**2.5.9. Corollary** Let $L$ be a free color Lie p-algebra freely generated by $X = \{x_1, \ldots, x_r\}$. Then

$$H(L, t) = -\sum_{n=1}^{\infty} \frac{\mu_p(n)}{n} \ln(1 - rt^n).$$

Suppose that $L$ is a free color Lie p-superalgebra generated by $X = \{x_1, \ldots, x_r\}$, and is multihomogeneous with respect to the set $X$. For elements of $L$ we introduce the multidegree $= (\alpha_1, \ldots, \alpha_r) \in \mathbb{N}_0^r$, and the degree $|\alpha| = \alpha_1 + \cdots + \alpha_r$. We have the following analogue of the Witt formula for the dimension of the multihomogeneous components of $L$.

**2.5.10. Corollary** *Let $L$ be a free color Lie p-algebra freely generated by $X = \{x_1, \ldots, x_r\}$. Then*

$$\dim L_n = \frac{1}{n} \sum_{m|\alpha} \mu_p(m) r^{\frac{n}{m}},$$

$$\dim L_\alpha = \frac{1}{|\alpha|} \sum_{m|\alpha} \mu_p(m) \frac{(|\alpha|/m)!}{(\alpha_1/m)! \ldots (\alpha_r/m)!}.$$

*When $L$ is the ordinary free Lie p-algebra, we get Petrogradsky's formulas [12].*

Petrogradsky initially proved the following Theorem for Lie superalgebras in [10, 11].

**2.5.11. Theorem** Let $L = L(X) = \bigoplus_{n=1}^{\infty} L_n$ be a free color Lie p-algebra ($G = G_+$) generated by a $\Lambda$-graded set $X = \bigcup_{\lambda \in \Lambda} X_\lambda$. Then

$$ch_\Lambda L_n = \frac{1}{n} \sum_{k|n} \mu_p(k) \left(ch_{[k]} X\right)_\Lambda^{\frac{n}{k}}.$$

*Proof.* We consider the new semigroup.

$$\Lambda' = \Lambda \times \mathbb{N}.$$

Define a weight function

$$wt : X \to \Lambda' : x \mapsto (\lambda, 1), x \in X_\lambda.$$

Then, we consider $L$ as a $\Lambda'$-graded. If we denote $t = e^{(0,1)}$ and $e^\lambda = e^{(\lambda,0)}$, then

$$\begin{aligned} ch_{\Lambda'} X &= \sum_{(\lambda, i) \in \Lambda'} |X_{(\lambda,i)}| e^{(\lambda,i)} \\ &= \sum_{\lambda \in \Lambda} |X_{(\lambda,1)}| e^{(\lambda,1)} \\ &= \sum_{\lambda \in \Lambda} |X_{(\lambda,1)}| e^{(\lambda,0)} e^{(0,1)} \\ &= t ch_\Lambda X. \end{aligned}$$

Using Theorem 2.5.7 and the operator of dilation, we see that.

$$\begin{aligned} ch_{\Lambda'} L &= -\sum \frac{\mu_p(k)}{k} \ln\left(1 - ch_\Lambda^{[k]} X\right)_{k=1}^\infty \\ &= -\sum \frac{\mu_p(k)}{k} \ln\left(1 - t^k ch_\Lambda^{[k]} X\right)_{k=1}^\infty. \end{aligned}$$

By the expansion of the logarithm, we have

$$ch_{\Lambda'} L = \sum_{k=1}^\infty \frac{\mu_p(k)}{k} \sum_{m=1}^\infty \frac{t^{mk}\left(ch_{[k]} X\right)^m}{m}.$$

Hence,

$$ch_{\Lambda'} L = \sum_{n=1}^\infty \frac{t^n}{n} \sum_{k|n} \mu_p(k) \left(ch_{[k]} X\right)^{\frac{n}{k}}.$$

However, it is undeniable that

$$ch_{\Lambda'}L = \sum_{n=1}^{\infty} ch_{\Lambda}L_n t^n.$$

Therefore,

$$ch_{\Lambda}L_n = \frac{1}{n}\sum_{k|n} \mu_p(k)\left(ch_{[k]}X\right)_{\Lambda}^{\frac{n}{k}},$$

as desired

Suppose that $G = G_+ = \{g_1, \ldots, g_r\}$ is an abelian group, and $L$ is a free color Lie p-superalgebra freely generated by a set $X = X_{g1} \cup \cdots \cup X_{gr}$ with $|X_{gi}| = s_i \geq 1\ i = 1, \ldots, r$. We define a weight function

$$wt: X \to \mathbb{N}^r : x \mapsto \lambda_i, for\ i = 1, \ldots, r\ and\ x \in X_{gi},$$

where $\lambda_i = (0, \ldots, 0, 1, 0, \ldots, 0)$ with 1 in the $i$th place. Again, we denote $t_i = e^{\lambda i}$, and so we have the following result.

**2.5.12. Theorem** *Suppose that $G = G_+ = \{g_1, \ldots, g_r\}$ is an abelian group, and $L$ is a free color Lie p-algebra freely generated by a set $X = X_{g1} \cup \cdots \cup X_{gr}$ with $|X_{gi}| = s_i \geq 1\ i = 1, \ldots, r$. Then*

1. $H(L; t_1, \ldots, t_r) = -\sum_{n=1}^{\infty} \frac{\mu_p(n)}{n} \ln\left(1 - \sum_{i=1}^{r} s_i t_i^n\right),$

2. $dimL_\alpha = \frac{1}{|\alpha|} \sum_{n|\alpha} \mu_p(n) \frac{\left(\frac{|\alpha|}{n}\right)!}{\left(\frac{\alpha_1}{n}\right)!\ldots\left(\frac{\alpha_r}{n}\right)!} s_1^{\frac{\alpha_1}{n}} \ldots s_r^{\frac{\alpha_r}{n}}.$

3. $dimL^{(n,g)} = \sum_{\substack{\alpha_1 + \cdots \alpha_r = n \\ g_1^{\alpha_1} \ldots g_r^{\alpha_r} = g}} dimL_{(\alpha_1, \ldots, \alpha_r)}.$

## 4. CONCLUSION

The Hilbert series of a Lie superalgebra representation is a generating function that encodes essential information about the representation. In particular, it provides a way to count the number of states in the representation with a given degree or energy. For common Lie superalgebra representations, Hilbert series has been computed explicitly in many cases, including for the basic representations of the Lie superalgebras. This series is a strong tool for comprehending the structure and properties of common Lie superalgebra representations. The computation of this topic relies on techniques from algebraic geometry and combinatorics, and it has important applications in mathematics.

## RECOMMENDATIONS

To strengthen the applications of Lie superalgebras. It will be helpful to give different methods for creating Hom-Lie superbialgebras. In addition, we could add research on triangular and coboundary Hom-Lie bialgebras. This algebra class is a generalization of both restricted Hom-Lie algebras and restricted Lie superalgebras [14]. In this paper, we could show how to get restricted Hom-Lie superalgebras from classical restricted Lie superalgebras using algebra endomorphisms.

## ACKNOWLEDGEMENTS

The author wishes to express his sincere gratitude to Dr. Yuri Bahturin, who initially suggested the subject matter. This research would not have been possible without his invaluable ideas. Additionally, the author extends his appreciation to his Ph.D. supervisor, Dr. Mikhail Kotchetov, for his invaluable guidance, insightful conversations, and helpful recommendations throughout the research process. Also, the author would like to thank Dr. Victor Petrogradsky for his suggestions. Their contributions have been instrumental in the completion of this work, and the author is deeply thankful for their support.